\documentclass{amsart}
\usepackage{amsfonts, amssymb}

\newtheorem{theorem}{Theorem}[section]

\newtheorem{corollary}[theorem]{Corollary}
\newtheorem{lemma}[theorem]{Lemma}
\newtheorem{proposition}[theorem]{Proposition}

\theoremstyle{definition}
\newtheorem{definition}[theorem]{Definition}
\theoremstyle{remark}
\newtheorem{remark}[theorem]{Remark}

\numberwithin{equation}{section}

\newcommand{\CC}{\mathbb C}

\newcommand{\aut}[1]{\mathrm{Aut\,}(#1)}
\newcommand{\dist}{\mathrm{dist }}
\newcommand{\supp}{\mathrm{supp }}

\def\cA{{\mathcal A}}

\def\cC{{\mathcal C}}

\def\cH{{\mathcal H}}

\def\cO{{\mathcal O}}

\def\re{\mathrm{Re}\,}
\def\dist{\mathrm{dist}\,}

\begin{document}

\title{Positivity and completeness of Invariant metrics}%

\author{Taeyong Ahn, Herv\'{e} Gaussier, Kang-Tae Kim}%
\address{(Ahn) Center for Geometry and its Applications, 
POSTECH, Pohang 790-784, The Republic of Korea}
\email{(Ahn) triumph@postech.ac.kr}
\address{\begin{tabular}{ll}
                     (Gaussier) & Univ. Grenoble Alpes, IF, F-38000 Grenoble, France\\
                     & CNRS, IF, F-38000 Grenoble, France
         \end{tabular}
}
\email{(Gaussier) herve.gaussier@ujf.grenoble.fr}
\address{(kim) Center for Geometry and its Applications 
and Department of Mathematics, POSTECH, Pohang 790-784, 
The Republic of Korea}
\email{(Kim) kimkt@postech.ac.kr}

\thanks{Research of the first and the third named authors is 
supported in part by the grant 2011-0030044 (The SRC-GAIA)
of the NRF of Korea.}%
\date{\today}

\begin{abstract}
We present a method for constructing global holomorphic peak 
functions from local holomorphic support functions for broad 
classes of unbounded domains in $\CC^n$. As an application, 
we establish a method for showing the positivity and completeness 
of invariant metrics including the Bergman metric mainly for the 
unbounded domains.
\end{abstract}
\maketitle

\section{Introduction}
The purpose of this article is twofold: (1) to present a method 
of obtaining global holomorphic peak functions at boundary points 
for unbounded domains in $\CC^n$ from any local holomorphic support 
functions, and (2) to establish the positivity and completeness of 
invariant metrics, primarily
of the Bergman metric, of certain unbounded domains in $\CC^n$.

As a consequence, we shall demonstrate applications of these methods 
to broad collections of unbounded domains that include 
the Kohn-Nirenberg domains, the Forn{\ae}ss domains, more 
generally those defined by (the ``positive variations'' of) 
weighted-homogeneous plurisubharmonic polynomial 
defining functions and some more. Notice that
several of these domains are not known whether they can be 
biholomorphic to bounded domains.

\section{Preliminary notation and terminology}

\subsection{The ball}
With $\|(z_1, \ldots, z_n)\|= \sqrt{|z_1|^2 + \ldots + |z_n|^2}$ 
we let 
$B^n(p,r) = \{z \in \CC^n \colon \|z-p\|<r\}$.   If $p$ is the origin 
and $r=1$, then we denote by $B^n := B^n (0,1)$.

\subsection{Distance to the complement}
$\delta_U (z) := \min \{1, \dist(z, \CC^n 
\setminus U)\}$, where $U$ is an open subset in $\CC^n$ and 
``dist'' means the Euclidean distance.

\subsection{Holomorphic peak and support functions}
Let $\Omega$ be an open set in $\CC^n$. Let 
$\cO(\Omega):=\{h\colon\Omega \to \CC,  h\colon 
\textrm{holomorphic}\}$ and $\cO^*(\Omega)
:=\{h\colon\Omega \to \CC, h\neq 0 \textrm{ on }\Omega, 
h\colon \textrm{holomorphic}\}$.
 
For a boundary point $p \in \partial\Omega$, a {\it peak function at 
$p$ for $\cO(\Omega)$} (or, a {\it global holomorphic peak function for 
$\Omega$ at $p$}) is defined to be a holomorphic function 
$f \in \cO(\Omega)$ such that:
\begin{itemize}
\item[(\romannumeral 1)] $\lim_{\Omega \ni z \to p} 
f(z) = 1$, and
\item[(\romannumeral 2)]  for every $r>0$ there exists $s>0$ 
satisfying $|f(z)|<1-s$ for every $z \in \Omega \setminus 
B^n (p,r)$. 
\end{itemize}
The point $p$ is called a {\it peak point} of $\Omega$ for $\cO(\Omega)$ 
(or, a {\it global holomorphic peak point of $\Omega$}) in such a case. 

By a {\it local holomorphic peak function} at a boundary point, say 
$q \in \partial\Omega$, we mean a peak function 
at $q$ for $\cO(V \cap \Omega)$ for some open neighborhood $V$ of $q$ 
in $\CC^n$. The point $q$ is then called a {\it local holomorphic peak 
point of $\Omega$}. 

If a holomorphic function $f$ is defined in an open neighborhood 
of the closure $\overline{\Omega}$ of the open set $\Omega$ in 
$\CC^n$ in such a way that it is also a peak function at $p$ for 
$\cO(\Omega)$, such $f$ is called a {\it global holomorphic support 
function} of $\Omega$ at $p$.  Local holomorphic support functions are 
defined likewise: a {\it local holomorphic support function} of $\Omega$ at 
$p \in \partial\Omega$ is a holomorphic function $g$ defined in an open 
neighborhood, say $V$, of $p$ in $\CC^n$ such that $g$ is a global 
holomorphic support function of $V \cap \Omega$ at $p$.

\section{Technical theorem for unbounded domains} 

Here we present the main technical theorem on how to obtain 
a global holomorphic peak function from a local holomorphic support 
function.

If $g\in \cO(V)$,
then we denote by
$
Z_g := \{ z \in V \colon g(z)=1\}
$.

\begin{theorem}\label{main_tech}
Let $\Omega$ be a domain in $\CC^n$.
If $p \in \partial\Omega$ satisfies the following two properties:
\begin{enumerate}
\item There exists an open neighborhood $V$ of $p$ in $\CC^n$ and a 
function
$g \in \cO(V)$ supporting $V \cap \Omega$ at $p$.
\item There are constants $r_1, r_2, r_3$ with $0<r_1<r_2<r_3<1$ and 
$B^n (p,r_3) \subset V$, and 
there exists a Stein neighborhood $U$ of $\overline\Omega$
and a function $h \in \cO(\Omega \cup V) \cap \cO^*(V)$ satisfying
$$
Z_g \cap U \cap \left(\overline{B^n (p,r_2)} \setminus B^n (p,r_1)\right) = \emptyset  
\eqno{(\dag)}
$$
and
$$
|h(z)|^2 \le C_0 \frac{\delta_U (z)^{2n}}{(1+\|z\|^2)^2},
\forall z \in \Omega,  \eqno{(\ddag)}
$$
\end{enumerate}
for some positive constant $C_0$, then $\Omega$ admits a peak function at $p$ for $\cO(\Omega)$.
\end{theorem}

\noindent\it Proof.
\rm
Take a $\cC^\infty$ function $\chi:\CC^n\to[0,1]$ satisfying 
$\chi\equiv 1$ on $B^n(p, r_1)$ and $\supp \chi\subseteq B^n(p, r_2)$.

Define a smooth $(0,1)$-from $\alpha$ on $U$ as follows 
(cf.\ \cite{FM}):
$$
\alpha (z) = \begin{cases} 
\bar\partial\Big(\frac{\chi(z)}{h(z)(1-g(z))} \Big)
& \textrm{if } z \in U \cap (B^n(p,r_2) \setminus 
\overline{B^n(p,r_1)}) 
\\
0 & \textrm{if } z \textrm{ is elsewhere in } U.
\end{cases} 
$$
By a theorem of H$\ddot{\rm o}$rmander (\cite{Hor_73}, 
Theorem 4.4.2), there exists a function $u \colon U \to \CC$ such 
that $\overline\partial u =\alpha$ on $U$ satisfying
$$
\int_{U} \frac{|u(z)|^2}{(1+\|z\|^2)^2} 
d\mu(z)
\le
\int_{U} |\alpha(z)|^2 d\mu(z),
$$
where $\mu$ denotes the Lebesgue measure for $\CC^n$. Notice that
$u \in C^\infty(U)$ by elliptic regularity. 

Since $\alpha$ is a bounded-valued smooth $(0,1)$-form with 
bounded support in $U$, the right-hand side of the preceding 
inequality is bounded above by a positive constant $C_1$, for 
instance. 
\medskip

Now we wish to obtain a pointwise estimate for $|u(z)|$.
\smallskip

Let $R>1$ be a constant such that $B^n(p, r_2)\subset B^n(0, R-1)$.
Let $\xi \in \Omega \setminus B^n (0,R)$. Our current aim is to
estimate $|u(\xi)|$.
\smallskip
 
For such $\xi$, we see that 
$\supp\chi \cap B^n(\xi, \delta_{U}(\xi)) = \emptyset$. This
implies that $u$ is holomorphic on $B(\xi, \delta_{U}(\xi))$.
These discussions yield 
\begin{eqnarray*}
C_1 & \ge & \int_{U} |\alpha(z)|^2 d\mu(z) 
\ge \int_{U} \frac{|u(z)|^2}{(1+\|z\|^2)^2} d\mu(z)
\\
& \ge & 
\int_{B(\xi, \delta_{U}(\xi))} 
\frac{|u(z)|^2}{(1+\|z\|^2)^2} d\mu(z)
\\
& \ge & 
\frac1{9(1+\|\xi\|^2)^2}\int_{B(\xi, \delta_{U}(\xi))} |u(z)|^2 d\mu(z)
\\
& \ge &
\frac{|u(\xi)|^2\cdot\mathrm{Vol}(B(\xi, \delta_{U}(\xi)))}{9(1+\|\xi\|^2)^2},
\end{eqnarray*}
where the last inequality is due to the sub mean-value
inequality.  In short, there exists a constant $C_2>0$ such that
$$
\Big(\frac{\delta_U (\xi)^{n}}{1+\|\xi\|^2}\Big)^2 |u(\xi)|^2 
\le C_2.
$$

At this stage we use the assumption (\ddag) and arrive at
$$
|h(\xi)u(\xi)|^2\leq C_3 \textrm{ for any } \xi \in \Omega 
\setminus B^n (0,R) 
$$
for some constant $C_3 > 0$.
\smallskip

Consider the case where $\xi \in \Omega \cap B^n (0,R)$.  
Notice that the assumption (\ddag) implies that $h$ is bounded. 
Since $u$ is smooth in $U$ and hence smooth at every point of 
$\overline\Omega$, the function $h(z)u(z)$ is bounded on 
$\Omega \cap B^n (0,R)$. Altogether, $hu$ is bounded on $\Omega$. 
\medskip

We are now going to construct a global holomorphic peak function
for $\Omega$ at $p$.  By the preceding arguments, we may choose 
a positive constant $c$ such that
$c |h(z) u(z)| < \frac12$ for any $z \in \Omega$.  
In particular, we obtain
$$
\re \big(c h(z) u (z) - 1\big) < -\frac12, 
\quad \forall z \in \Omega.
$$
Consider 
$$
\psi (z) = \begin{cases} 
\frac{\chi(z)}{1-g(z)}
& \textrm{if } z \in \Omega \cap B^n(p,r_2) 
\\
0 & \textrm{if } z \in \Omega \setminus B^n(p,r_2).
\end{cases} 
$$
Note that $\psi(z)-h(z)u(z)$ is a holomorphic function on 
$\Omega$ since $\overline\partial u=\alpha$. Altogether, 
if we define $f\colon\Omega\to\CC$ by
$$
f(z) := \exp\left((-c[\psi(z)-h(z)u(z)]-1)^{-1}\right),
$$
then $f$ is holomorphic on $\Omega$.
\medskip

Note that the real part of $c h(z)u (z)-1$ on $\Omega$ is 
negative. The function $-\re c \psi$ on $\Omega$ is also negative
due to its construction. Hence
$$
\re (-c[\psi(z)-h(z)u(z)]-1)^{-1}<0.
$$
So $|f(z)|<1$ for any $z \in \Omega$.

Finally, $\lim_{\Omega \ni z\to p} f(z) = e^0 = 1$ since 
$\psi$ tends to $\infty$ as $z \to p$. Notice that $p$ is the 
only boundary point that has this property for $f$ by the 
construction throughout. 

The remaining condition for $f$ to be a peak function at $p$ for 
$\cO(\Omega)$ is also easily checked from the definition of $f$ itself. 
We now see that $f$ is the desired global holomorphic peak function for 
$\Omega$ at $p$, and hence the proof is complete. \hfill $\Box$
\bigskip

\begin{remark} \rm 
Note that condition ($\ddag$) is necessary only for the case where 
$\Omega$ is unbounded.  If $\Omega$ is bounded, the assumption 
($\ddag$) is void; one may just take $h\equiv 1$.
\end{remark}

\section{Applications to unbounded domains}

\subsection{Weighted homogeneous domains in $\CC^n$}
For $z=(z_1, \ldots, z_n) \in \CC^n$, denote by 
$z'= (z_1, \ldots, z_{n-1})$ and consequently $z=(z',z_n)$.

Let $m_1, ..., m_{n-1}$ be positive integers. A real-valued polynomial 
$P$ on $\CC^{n-1}$ is called 
weighted-homogeneous of weight $(m_1, ..., m_{n-1})$, if 
$$
P\big(t^{\frac1{2m_1}} z_1, ..., t^{\frac1{2m_{n-1}}} z_{n-1}\big)
= tP(z_1, ..., z_{n-1}),\ \forall t>0,\ \forall(z_1,\dots,z_{n-1}) \in \CC^{n-1}.
$$
If $m=m_1=\ldots=m_{n-1}$ then $P$ is called {\it homogeneous of 
degree $m$}. 

\begin{definition}\label{def:WP}
A domain $\Omega$ in $\CC^n$ is called a {\it WB-domain} 
(meaning ``weighted-bumped'') if 
$$
\Omega = \{ z \in \CC^n \colon \re z_n + P(z') < 0 \},
$$
where: 
\begin{itemize}
\item[(\romannumeral 1)] $P$ is a real-valued, weighted-homogeneous 
polynomial on $\CC^{n-1}$ of weight 
$(m_1, ..., m_{n-1})$,   
\item[(\romannumeral 2)] $P$ is plurisubharmonic (psh) without 
pluriharmonic terms, and
\item[(\romannumeral 3)] there is a constant $s > 0$ such that 
$P(z')-2s\sum_{j=1}^{n-1} |z_j|^{2m_j}$ is also psh in $\CC^{n-1}$.
\end{itemize}
\end{definition}

We prove the following theorem.

\begin{theorem}\label{wb-homo}
If $\Omega$ is a WB-domain in $\CC^n$, there exists  
a peak function at $0=(0,\ldots,0)$ for $\cO(\Omega)$,
continuous on $\overline\Omega$ with exponential decay
at infinity and nowhere zero.
\end{theorem}

\noindent
\it
Proof.
\rm
We shall use the notation of Definition \ref{def:WP}.
Let $\displaystyle \mu_j := \frac1{m_j} \prod_{k=1}^{n-1}m_k$ for 
$j=1, \ldots, n-1$ and let $H(z_1, \cdots, z_{n-1}):=P(z_1^{\mu_1}, \cdots, 
z_{n-1}^{\mu_{n-1}})$. Then $H$ is a homogeneous polynomial of 
degree $2k = 2m_1\cdots m_{n-1}$.
\medskip

Let $\Omega_H:=\{(z_1, \cdots, z_n):\re z_n+H(z_1, \cdots, 
z_{n-1})<0\}$ and let $F_H(z_1, \cdots, z_n) = 
(z_1^{\mu_1}, \cdots, z_{n-1}^{\mu_{n-1}}, z_n)$. 
Then, $F_H\colon \Omega_H \to \Omega$ is a holomorphic 
ramified covering map of finite degree.

By Theorem 4.1 of \cite{Bedford-Fornaess}, the domain 
$\widetilde\Omega:=\{(z_1, \cdots, z_n)\colon \re z_n+\widetilde 
H(z_1, \cdots, z_{n-1}) -\delta|z_n| - \delta \sum_{j=1}^{n-1} |z_j|^{2k}<0\}$,
for some $\delta > 0$,  
admits a peak function, which we denote by $Q_H(z)$ here, at the 
origin for $\cO(\widetilde\Omega)$.  This peak function by Bedford 
and Forn{\ae}ss enjoys an exponential decay condition at infinity and 
vanishes nowhere on $\widetilde\Omega$.

To obtain the desired peak function at $0$ for $\cO(\Omega)$, we 
symmetrize $Q_H$: 
for each $z=(z_1, \cdots, z_n)\in\overline\Omega$, there are 
exactly $\prod_{j=1}^{n-1}\mu_j$ preimages by the map $F_H$
(counting with multiplicity), of the following form:
$$
\big( |z_1|^{\frac1{\mu_1}}
e^{\frac{\sqrt{-1}}{\mu_1}{(2k_1\pi+\arg z_1)}}, 
\ldots, |z_{n-1}|^{\frac1{\mu_{n-1}}} 
e^{\frac{\sqrt{-1}}{\mu_{n-1}} (2k_{n-1}\pi+\arg z_{n-1})}, z_n 
\big),
$$
where $k_j \in \{1, \ldots, \mu_j\}$ for every $j$.
Let  
$$
q(w):=
\prod_{\substack{1\le j\le n-1 \\ 1\le k_j\le \mu_j}} 
Q_H\big(e^{\frac{2k_1 \sqrt{-1} \pi}{\mu_1}} w_1, \cdots, 
e^{\frac{2k_{n-1} \sqrt{-1} \pi}{\mu_{n-1}}} w_{n-1}, w_n \big).
$$
Then $q\colon \Omega_H \to \CC$ is a well-defined holomorphic 
function on $\Omega_H$. Moreover, this defines a unique
holomorphic function $Q \in \cO(\Omega) \cap \cC^0 (\overline\Omega)$ 
satisfying $Q\circ F_H (w) = q(w)$ for any $w \in Q_H$. 
It is now immediate that this $Q$ is the desired peak function at 
$0$ for $\cO(\Omega)$.
\hfill $\Box$
\bigskip

More importantly for our purpose, this peak 
function $Q$ is holomorphic on the domain
$$
\Omega^\epsilon := \{z \in \CC^n \colon \re z_n + P(z') < \epsilon|z_n|+\epsilon 
\sum_{j=1}^{n-1} |z_j|^{2m_j} \}  
$$
and has an exponential decay condition at infinity.  It also 
vanishes nowhere.  Hence it can play the role of the function
$h$ of Theorem \ref{main_tech}.
As a consequence we obtain (using the homothety and translation 
automorphisms for weakly pseudoconvex points, if necessary) 
the following

\begin{corollary} \label{WB-peak}
If $\Omega$ is a WB-domain then every boundary point admits a peak
function for $\cO(\Omega)$.
\end{corollary}

\begin{remark} \rm
Note that this corollary generalizes the case of domains 
with defining functions of diagonal type studied by Herbort in
\cite{Herbort}. The method of this article yields an alternative 
proof to Theorem 2 of \cite{Herbort}.

\end{remark} 

\begin{remark}
Theorem \ref{wb-homo} can be understood as a generalization of 
Theorem 4.1 of \cite{Bedford-Fornaess}. There is another 
generalization of it in a different direction by Noell (\cite{Noell}). 
For $P$ with the assumptions by Noell in \cite{Noell}, the above 
arguments still hold and we can obtain a statement corresponding 
to Corollary \ref{WB-peak}.

\end{remark}

\subsection{The Kohn-Nirenberg domains, the Forn{\ae}ss domains 
and positive variations of WB-domains} \label{KN}
The domains
$$
\Omega_\textrm{HKN} = \{(z,w) \in \CC^2 \colon \re w + |z|^8 
+ \frac{15}7 |z|^2 \re z^6 < 0 \}.
$$
and
$$
\Omega_\textrm{KN} = \{(z,w) \in \CC^2 \colon \re w + |zw|^2
+ |z|^8 +  \frac{15}7 |z|^2 \re z^6 < 0 \}.
$$
were first introduced in \cite{KN_73}; for $\Omega_\textrm{HKN}$ 
and $\Omega_\textrm{KN}$, the origin is 
the boundary point that does not admit, even locally, any 
holomorphic support functions, despite the fact that the boundary 
is real-analytic everywhere and strongly pseudoconvex at every
boundary point except the real line $\{z=0, \re w=0\}$ for 
$\Omega_\textrm{HKN}$, and except the origin for 
$\Omega_\textrm{KN}$, respectively. There are still some problems yet 
to be answered for these domains, as recent research concerning 
unbounded domains attracts much attention (\cite{CKO}, \cite{HST} 
and also \cite{Herbort}, \cite{Herbort1}).

Note that $\Omega_\textrm{HKN}$ belongs to the class of WB-domains. 
Hence, we will focus on $\Omega_\textrm{KN}$ and call it the Kohn-
Nirenberg domain in the rest of the article.

\begin{theorem}\label{peak-thm}
There is a peak function for $\cO(\Omega_\mathrm{KN})$ 
at every boundary point of $\Omega_\mathrm{KN}$.
\end{theorem}

\noindent
\it
Proof.
\rm
Observe that $\Omega_\textrm{KN}\subset \Omega_\textrm{HKN}$ and
that 
\( 0 \in \partial \Omega_\textrm{KN}\cap 
\partial \Omega_\textrm{HKN} \). 
Consider now the domain $\Omega_\textrm{HKN}$.  Then 
Theorem \ref{wb-homo} provides a special peak
function at the origin for this domain.  It continues to be a peak 
function for $\cO(\Omega_\textrm{KN})$ at the origin.  This peak 
function also plays the role of $h$ in the hypothesis of 
Theorem \ref{main_tech}. For $U$, we simply take
$$
U=\{(z,w) \in \CC^2 \colon \re w + |zw|^2
+ |z|^8 +  \frac{15}7 |z|^2 \re z^6 < \varepsilon \}
$$
where $\varepsilon$ is chosen small enough.
Then, since the defining function of $U$ is a polynomial and $h$ 
decays exponentially, the conditions of Theorem \ref{main_tech} are 
satisfied. Recall that every boundary point except the origin is a 
strongly pseudoconvex point. Hence the assertion follows immediately. 
\hfill $\Box$
\medskip

We remark that the 
same applies to the Forn{\ae}ss domains in $\CC^2$ 
(\cite{Fornaess_77}) defined by
$$
\Omega_\textrm{HF} := \{(z, w)\in\CC^2\colon\re w + |z|^6 
+ t |z|^2 \re z^4 < 0\}
$$
and
$$
\Omega_\textrm{F}:=\{(z, w)\in\CC^2 \colon
\re w + |zw|^2 + |z|^6 + t |z|^2 \re z^4 < 0\}
$$
for constant $t$ with $1<t<\frac95$.
\medskip

Since the arguments of this type have a general nature, we formulate
it into the following formal statement:  

\begin{theorem} \label{W_S}
Let $W := \{z \in \CC^n \colon \re z_n + P(z')<0\}$ be a WB-domain. 
Let $S$ be a non-negative plurisubharmonic polynomial defined on 
$\CC^n$ such that $S(0)=0$ and 
$$
W_S := \{z \in \CC^n \colon 
\re z_n + S(z) + P (z') < 0\}.
$$
If a boundary point of $W_S$ admits a local holomorphic support 
function, then it also admits a global holomorphic peak function.
\end{theorem}

\begin{remark}\label{rem:EXP} 
Using this type of argument, we can also handle certain unbounded 
domains which are subdomains of WB-domains. In particular, we can 
handle some unbounded domains defined by a non-polynomial defining 
function. For example, consider
$$
W_E:=\{(z, w)\in\CC^2: \re w+ \exp(|z|^2)<0\}.
$$
Note that this domain has infinite volume. As the open set $U$ of 
Theorem \ref{main_tech}, we take
$$
U:=\{\re w+ \exp(|z|^2)-\varepsilon \exp(|z|^2)< \varepsilon\}
$$
with $\varepsilon$ small enough.

Moreover, using the inequality $|z|^2<\exp(|z|^2)$, we can take 
as the function $h$ of Theorem \ref{main_tech} the one given by 
Theorem \ref{wb-homo}. Since every boundary point of $W_E$ 
is a strongly pseudoconvex point, every boundary point of $W_E$ 
admits a global holomorphic peak function. Indeed, in this case, 
they are all global holomorphic support functions. Another 
example with different nature will be dealt with in the next subsection.
\end{remark}

\subsection{Some unbounded domains with finite volume} The domain
$$
E := \{(z, w) \in \CC^2 \colon |w|^2<\exp(-|z|^2)\}
$$
cannot be biholomorphic to any bounded domain since it contains
the complex line $\{(z,w) \in \CC^2 \colon w=0\}$.  However, it is
strongly pseudoconvex at every boundary point.  

We point out that Theorem \ref{main_tech} applies to this 
domain as well by setting for $\varepsilon$ small enough
$$
U := \{(z, w) \in \CC^2 \colon |w|^2<\exp(-|z|^2)+\varepsilon\}
$$
and $h(z,w):=w$.
Consequently, every boundary point admits a peak function for 
$\cO (E)$, which is, especially in this case, a global holomorphic 
support function.

\section{Applications to invariant metrics of unbounded domains}

\subsection{Remarks on the Hahn-Lu comparison theorem}

Recall the following classical concepts: for a complex 
manifold $M$, denote by $\cO(M, B^1)$ the set of 
holomorphic functions from $M$ into the unit open disc 
$B^1$.  Let $p \in M$ and $v \in T_p M$. Then the 
{\it Caratheodory pseudo-metric} ({\it metric}, if positive) 
of $M$ is defined 
by
$$
{\frak c}_M (p,v) 
= \sup \{ |df_p (v)| \colon f \in \cO(M, B^1), f(p)=0 \}.
$$
This induces the {\it integrated Caratheodory pseudo-distance} 
({\it distance}, if positive) 
$$
\rho^{\frak c}_M (p, q) = \inf \int_0^1 {\frak c}_M 
(\gamma(t), \gamma'(t)) dt,
$$
where the infimum is taken over all the piecewise $\cC^1$ curves
$\gamma\colon[0,1] \to M$ with $\gamma (0)=p$, $\gamma (1)=q$.

For the Bergman metric, one starts with the space 
$\cA^2 (\Omega) := \{ f \in \cO(\Omega) \colon 
\int_\Omega |f|^2 < +\infty \}$.  Equipped with the 
$L^2$-inner product, it is a separable Hilbert space. If one 
considers a complete orthonormal system $\{\varphi_j \colon j=1,2,
\ldots\}$ of $\cA^2(\Omega)$, the Bergman kernel function can be 
expressed as 
$K_\Omega (z,w) = \sum_{j \geq 1} \varphi_j(z)\overline{\varphi_j (w)}$. 
Then it defines 
$$
{\frak b}^\Omega_{ab} (\zeta) = \frac{\partial^2}{\partial z_a 
\partial \bar z_b} \Big|_\zeta \log K_\Omega (z, z)
$$
as well as the (1,1)-tensor
$$
{\frak b}^\Omega_\zeta := \sum_{a,b=1}^n {\frak b}^\Omega_{ab} 
(\zeta)\ d\zeta_a \otimes d\bar\zeta_b.
$$
If $K_\Omega (z, z)$ is non-zero, than ${\frak b}^\Omega_z$ 
defines a smooth (1,1)-Hermitian form that is positive 
semi-definite.  It is a result of Bergman himself that 
${\frak b}^\Omega$ is a positive definite Hermitian metric for 
bounded domains.  
In general this may not even be defined, and even when it is 
defined, it may not be positive.

Nevertheless, we shall follow the convention and write as 
${\frak b}^\Omega$ the {\it Bergman metric} of the domain 
$\Omega$ and the notation ${\frak b}^\Omega_p$ shall always
mean the Bergman metric of the domain $\Omega$ at the point 
$p\in\Omega$.
\medskip

We now present a modification of the comparison
theorem by Hahn \cite{Hahn1, Hahn_77, Hahn2}, and Lu \cite{Lu} which
compares the Caratheodory metric and the Bergman metric (even if 
both may be degenerate). 

\begin{theorem}[The Hahn-Lu comparison theorem] \label{HL}
If $M$ is a complex manifold and $p$ is a point in $M$ 
such that its Bergman kernel $K_M$ satisfies $K_M (p, p) \neq 0$,
then its Bergman metric ${\frak b}^M_p (v, w)$ and the Caratheodory 
pseudometric ${\frak c}_M (p, v)$ satisfy the inequality
$$
\big({\frak c}_M (p, v)\big)^2 \le {\frak b}^M_p (v, v), 
$$
for any $v \in T_p M$. 
\end{theorem}

\noindent\it Proof. \rm We shall only prove it for the 
case when $M=\Omega$
is a domain in $\CC^n$, staying closely to the purpose 
of this article; the manifold case uses essentially the 
same arguments except some simplistic adjustments. 
  
Start with the following quantities developed
by Bergman \cite{Bergman_35}: 
\begin{eqnarray*}
{\frak B}_0 (p) 
& = & \sup\Big\{|\psi (p)|^2 \colon \psi 
\hbox{ holomorphic}, \int_\Omega |\psi|^2 \le 1 \Big\}
\\
{\frak B}_1 (p, v) 
& = & \sup\Big\{ \Big|\partial_v \varphi |_{p}\Big|^2 
\colon \varphi \textrm{ holomorphic}, \varphi(p)=0,
\int_\Omega |\varphi|^2 \le 1 \Big\},
\end{eqnarray*}
where $\partial_v \varphi|_p = \sum_{j=1}^n v_j 
\frac{\partial \varphi}{\partial z_j}\Big|_{p}$. 
These concepts are significant because 
${\frak B}_1 (p, v) = {\frak B}_0 (p) \cdot
{\frak b}^\Omega_p (v, v)$, when ${\frak B}_0 (p)>0$.

By the Cauchy estimates and Montel's theorem, there exists 
an $L^2$-holomorphic function $\hat\psi$  on $\Omega$ with 
$\|\hat\psi\|_{L^2(\Omega)} \le 1$ satisfying
$|\hat\psi (p)|^2 = \frak{B}_0(p)$. (See \cite{Hahn_77, GKK}).
Then Montel's theorem on normal families implies the existence of
$\eta \in \cH(\Omega, B^1)$ on $\Omega$ with 
$\eta(p)=0$ and $\big|\partial_v \eta|_{p}\big|=
|d\eta_{p} (v)| = {\frak c}_\Omega (p, v)$, the 
Caratheodory length of $v$ at $p$.  Since
$|\eta\hat\psi| \le |\hat\psi|$, $\|\eta\hat\psi\|_{L^2(\Omega)} 
\le \|\hat\psi\|_{L^2(\Omega)} \le 1$. Since $\eta(p)=0$, 
$\big|\partial_v  (\eta\hat\psi)|_{p}\big|
= \big|\partial_v  \eta|_{p}\big|~ |\hat\psi(p)|$.  Altogether, we 
arrive at
$$
\frak{B}_1 (p,v) \ge |\partial_v (\eta \hat\psi)|_p |^2
= {\frak c}_\Omega (p, v)^2 ~\frak{B}_0(p).
$$
This immediately yields the comparison inequality
$$
{\frak b}^{\Omega}_p (v, v) = \frac{\frak{B}_1 (p,v)}{\frak{B}_0(p)} \ge 
{\frak c}_\Omega (p, v)^2.
$$
\hfill $\Box$
\medskip

\begin{remark} \rm
The original statements required positivity of both 
metrics. But the proof above (almost identical with the 
arguments by Hahn \cite{Hahn2}) clearly shows that not all those 
assumptions are necessary. On the other hand, this modification is 
significant since one obtains that the Bergman metric is positive
definite whenever the manifold is Caratheodory hyperbolic.  We shall 
see applications in the next section.
\end{remark}

\subsection{Positivity and completeness of invariant metrics}

\subsubsection{Positivity}
As a result of discussions above, we present the following:

\begin{proposition}
The Kohn-Nirenberg domains, the Forn{\ae}ss domains, all
WB-domains as well as the domain $W_S$ in Theorem \ref{W_S} 
are Caratheodory hyperbolic.  Consequently, their Kobayashi 
metric and Bergman metric are positive. 
\end{proposition}

\noindent\it Proof. \rm  
Let $Q$ be an appropriate global holomorphic peak function at 
the origin for each case introduced in Section 4. 
For $p = (p_1, p_2)$ and $v=(v_1,v_2)$, take $g(z_1, z_2) := 
Q(z_1, z_2) \cdot \big(\bar v_1 (z_1-p_1)
+ \bar v_2 (z_2-p_2) \big)$.  Then $g$ is a bounded holomorphic 
function on $\Omega$, since $Q$ decays exponentially at infinity in 
$\Omega$ and is continuous in the closure of $\Omega$. Moreover 
$g(p) = 0$ and $|dg_p (v)| = Q(p)\|v\|^2 > 0$.  Hence the 
Caratheodory metric is positive on $\Omega$.  Consequently the Kobayashi metric, 
being larger, is positive as well on $\Omega$. For the Bergman metric, we use the 
fact that $Q$ belongs to $\cA^2(\Omega)$ thanks to its exponential 
decay (whereas the domain is defined by a polynomial inequality). 
Thus the diagonal of the Bergman kernel vanishes nowhere on $\Omega$ and it follows from the Hahn-Lu comparison theorem (Theorem~\ref{HL}) that the Bergman metric is positive definite on $\Omega$.
\hfill $\Box$ 

\subsubsection{Completeness} Recall that a metric is said to 
be {\it complete} if every Cauchy sequence converges with respect to 
its integrated distance. Note that it is well-known that the Caratheodory metric is continuous with respect to the Euclidean metric.
Since the domains we handle are unbounded, we need to handle 
points near infinity in each of the following cases:
\bigskip

\noindent
{\it Case 1. WB-domains}: In this case we first prove

\begin{lemma}\label{infinity}
If $\Omega$ is a WB-domain equipped with its integrated 
Caratheodory distance $\rho^{\frak c}_\Omega$, then every 
Cauchy sequence is bounded away from infinity.
\end{lemma}

\noindent\it Proof. \rm 
Suppose that $(q_\nu)$ be a Cauchy sequence. Passing 
to a subsequence, without loss of generality, we may assume that 
$\lim_{\nu \to \infty} \|q_\nu\| = \infty$ to get a contradiction. 
Since WB-domains admit homothety automorphisms shrinking 
to the origin, one can 
always find $\varphi_\nu \in \aut \Omega$ such that $\|\varphi_\nu 
(q_\nu)\| = 1$ for every $\nu$.  Of course we have $\lim_{\nu \to 
\infty} \|\varphi_\nu (q_1)\| = 0$. Since the origin has a peak 
function, say $f$, for $\cO(\Omega)$ (cf.\ Theorem \ref{wb-homo}), 
it follows that
\begin{align*}
\lim_{\nu \to \infty} \rho^{\frak c}_\Omega (q_1, q_\nu)
& = \lim_{\nu \to \infty} \rho^{\frak c}_\Omega (\varphi_\nu(q_1), 
\varphi_\nu(q_\nu)) 
\\
& \ge \lim_{\nu \to \infty} d^P (f\circ\varphi_\nu(q_1),
f\circ\varphi_\nu(q_\nu)) = \infty,
\end{align*}
where $d^p(\cdot, \cdot)$ denotes the Poincar\'e distance of $B^1$. 
Of course this is impossible for a $\rho^{\frak c}_\Omega$-Cauchy 
sequence, and hence the proof is complete.
\hfill
$\Box$
\bigskip

Now, Corollary \ref{WB-peak}, Theorem \ref{HL} and 
Lemma \ref{infinity} imply

\begin{theorem}
If $\Omega$ is a WB-domain, then its Caratheodory metric is complete.  
Moreover, its Kobayashi metric and Bergman metric are also complete.
\end{theorem}
\bigskip

\noindent
{\it Case 2. The Kohn-Nirenberg domains, the Forn{\ae}ss domains and 
positive variations of WB-domains admittiong a local holomorphic support function at each boundary point}:  Let $\Omega$ denote any of these domains.
If $\Omega$ is either a Kohn-Nirenberg domain or a Fornaess domain, then by Theorem~\ref{peak-thm} the Caratheodory distance $\rho^{\frak c}_M(p,q_j)$ tends to infinity when $p \in \Omega$ and $(q_j)_j$ is a sequence of points in $\Omega$ converging to a point $q \in \partial \Omega$. If $\Omega$ is a positive variation of WB-domains admitting a local holomorphic support function at each boundary point, the same holds according to Theorem~\ref{W_S}.
Thus we only need to prove in each case that every Cauchy sequence $(z^j)_j$ of points in $\Omega$, with respect to the integrated Caratheodory distance $\rho^{\frak c}_{\Omega}$, is bounded away from infinity. Since $\Omega$ is a subdomain of a WB-domain $\Omega'$, it follows from the distance decreasing property of the Caratheodory metric that $(z^j)_j$ is a Cauchy sequence with respect to $\rho^{\frak c}_{\Omega'}$. Then $(z^j)_j$ is bounded away from infinity by Lemma~\ref{infinity}.
 Hence $\Omega$ is complete with respect to the Caratheodory metric. 
Using the comparisons of 
invariant metrics including the Hahn-Lu comparison theorem, the 
Bergman metric and the Kobayashi metric of $\Omega$ are also positive and complete. \qed

\medskip
\begin{remark} It is not difficult to see that the Caratheodory metric, 
Kobayashi metric and Bergman metric of the domain $W_E$ are also 
positive and complete.
\end{remark}
\bigskip

\noindent
{\it Case 3. $\Omega = \{(z,w) \in \CC \colon |w|^2 < e^{-|z|^2} \}$}:
\rm 
This unbounded domain does not admit any translation or homothety 
automorphisms.  Also it is not biholomorphic to any bounded domains. 
Also the Kobayashi metric degenerates at points on the complex line 
defined by the equation $w=0$, which we denote by $L$. 

However a direct computation shows that every holomorphic polynomial 
on $\CC^2$ belong to $\cA^2(\Omega)$ and the Bergman kernel is well-defined and positive along the diagonal. Indeed, a direct computation yields the 
Bergman kernel for $\Omega$
$$
K_\Omega(z, w)
=\frac{2\exp(2z_1\overline w_1) (1 +z_2\overline w_2 \exp(2z_1
\overline w_1))}{\pi^2(1-z_2\overline w_2\exp(2z_1\overline w_1))^3}.
$$
Also it is not hard to see
that the Bergman metric is positive definite.

Although the completeness of the Bergman metric of $\Omega$ 
follows by an explicit computation, 
Theorem 3.1 and Theorem 5.1 yield a different proof 
of that result. Since every boundary point admits a global 
holomorphic peak function and our version of Hahn-Lu 
comparison theorem applies pointwise, it suffices to check 
whether every Cauchy sequence with respect to the Bergman 
metric is bounded away from infinity. This is true since from 
the direct computations we see that the length of any piecewise 
$\cC^1$-curve always has larger length with respect to the 
Bergman metric than its projected image onto the complex line 
$L$ (with respect to the natural projection). Moreover, again 
from direct compuations, the Bergman metric restricted to the 
complex line $L$ is Euclidean. So, this proves the completeness 
of the Bergman metric.

\end{document}